\documentclass[a4paper,10pt]{article}
\usepackage[footnotesize]{caption}
\usepackage[caption=false,font=footnotesize]{subfig}
\usepackage{mathptmx,graphicx,amsmath,amsfonts,amsthm,amssymb,url,cite}
\usepackage[margin=0.7in]{geometry}
\usepackage[all]{xy}

\newtheorem{theorem}{Theorem}[section]

\begin{document}

\title{The Collatz conjecture and De Bruijn graphs}
\author{Thijs Laarhoven\thanks{Department of Mathematics and Computer Science, Eindhoven University of Technology, P.O. Box 513, 5600 MB Eindhoven, The Netherlands. E-mail: $\{$t.m.m.laarhoven,b.m.m.d.weger$\}$@tue.nl} \and Benne de Weger\footnotemark[1]}
\date{\today}
\maketitle

\abstract{We study variants of the well-known Collatz graph, by considering the action of the $3n + 1$ function on congruence classes. For moduli equal to powers of $2$, these graphs are shown to be isomorphic to binary De Bruijn graphs. Unlike the Collatz graph, these graphs are very structured, and have several interesting properties. We then look at a natural generalization of these finite graphs to the $2$-adic integers, and show that the isomorphism between these infinite graphs is exactly the conjugacy map previously studied by Bernstein and Lagarias. Finally, we show that for generalizations of the $3n + 1$ function, we get similar relations with $2$-adic and $p$-adic De Bruijn graphs. }


\section{Introduction}

The $3n + 1$ or Collatz conjecture is a long-standing open problem in mathematics. Let the $3n + 1$ function $T$ be defined on the integers by
\begin{align}
T(n) = \begin{cases} (3n + 1)/2 & \text{if $n$ is odd,} \\ n/2 & \text{if $n$ is even}. \end{cases}
\end{align}
The Collatz conjecture states that, starting from any positive integer $n$, repeated application of the function $T$ will eventually produce the number $1$, after which it will end in the cycle $\{1,2\}$. This conjecture is true if and only if, on the positive integers, there are no divergent paths (i.e., $\lim_{k \to \infty} T^k(n) < \infty$ for all positive integers $n$, where $T^0(n) = n$ and $T^{k+1}(n) = T(T^k(n))$ for $k \geq 0$) and there are no other cycles besides the trivial cycle $\{1,2\}$ (i.e., there are no natural numbers $n \geq 3$ with $T^k(n) = n$ for some $k \geq 1$). Though easy to state, this problem seems very hard, if not impossible to solve.

Because of its simple formulation, researchers from many different branches of mathematics have at one time or another encountered this problem and have become fascinated by it. This has lead to hundreds of papers in the last few decades, with each researcher using his own area of expertise to shed a new light on this problem. An excellent overview of many of these papers was given by Lagarias \cite{lagarias11, lagarias12}, while extensive surveys of previous work on this problem can be found in books by Lagarias \cite{lagarias10} and Wirsching \cite{wirsching98}.

Three of those branches of mathematics that have been used to study the Collatz conjecture are those of graph theory, modular arithmetic and $2$-adic integers. This paper aims to show connections between these three approaches.

We start in Section~\ref{sec:finite} with studying modular Collatz graphs, i.e., finite graphs that capture the behaviour of the $3n + 1$ function on congruence classes of integers. It turns out that there is an intimate relation to binary De Bruijn graphs when the modulus is a power of $2$. Letting this modulus grow to infinity, in Section~\ref{sec:infinite} we are led to studying these problems on the $2$-adic integers. This leads to a natural generalization of the binary De Bruijn graphs to the $2$-adic integers. In Section~\ref{sec:struct} we look at the structure of this infinite graph, and we try to describe how various Collatz graphs are embedded in it. In Section~\ref{sec:generalizations} we briefly indicate possible generalizations.


\section{Binary Collatz graphs and binary De Bruijn graphs}
\label{sec:finite}

One particular approach to the $3n + 1$ problem that caught our attention is using directed graphs to visualize the action, and in particular iteration, of the function $T$. We denote a directed graph by $G = (V,E)$, where $V$ is the set of vertices, and $E$ is the set of directed edges. Since we will not be dealing with undirected graphs, we will refer to directed graphs simply as graphs. Consider the graph $C(\mathbb{N}_{+}) = (V,E)$ with vertices $V = \mathbb{N}_{+}$ and edges $E = \{n \to T(n): n \in V\}$. This graph is known in the literature as the Collatz graph \cite{andaloro02, lagarias85b, monks12, wirsching98}, and allows for a simple visual explanation of the $3n + 1$ problem to a broad audience \cite{munroe10}. 

\begin{figure}[t]
\centering
\begin{align*}
\xymatrix{
*+[F-]{0} \ar@(l,d)[] \ar[r] & *+[F-]{2} \ar@(ld,rd)[] \ar@/^/[r] & *+[F-]{1} \ar@/^/[l]
}
\end{align*}
\caption{The modular Collatz graph with modulus $m = 3$.\label{fig1}}
\end{figure}
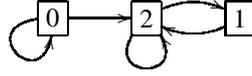

Since the Collatz graph has infinitely many vertices and looks very chaotic, we introduce a new family of related, finite graphs. Given some modulus $m$, we define the \textit{modular Collatz graph with modulus $m$} as the graph $G = (V, E)$ with vertices $V = \{0, 1, \ldots, m - 1\}$, and a directed edge runs from $a$ to $b$ if there exist some numbers $a_1, b_1 \in \mathbb{Z}$, with $a_1 \equiv a \pmod m$ and $b_1 \equiv b \pmod m$, such that $T(a_1) = b_1$. For instance, taking $m = 3$ leads to the graph on three vertices in Figure~\ref{fig1}.

As can be seen from this graph, once we `leave' the set $3\mathbb{Z}$ we will never return, as there are no edges going into the vertex $0$. The only way not to leave this set is to have $T^{k}(n) \equiv 0 \pmod 6$ for all $n$, which implies $n = 0$. Hence proving that all positive integers $n \equiv \pm 1 \pmod 3$ iterate to $1$ is sufficient to prove the Collatz conjecture. Although this conclusion may seem trivial, and similar but stronger results have been derived by the Monkses \cite{monks06, monks12}, it shows that studying these graphs may be useful.

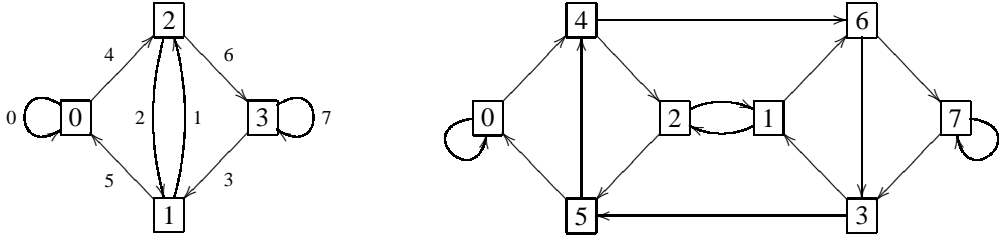
\begin{figure}[t]
\centering
\begin{align*}
\xymatrix{
 & *+[F-]{2} \ar[rd]^{6} \ar@/_/[dd]_{2} & \\
*+[F-]{0} \ar[ru]^{4} \ar@(ul,dl)[]_{0} &  & *+[F-]{3} \ar@(ur,dr)[]^{7} \ar[dl]^{3} \\
 & *+[F-]{1} \ar[ul]^{5} \ar@/_/[uu]_{1} &
} \quad \quad \quad \quad \quad \xymatrix{
 & *+[F-]{4} \ar[rrr] \ar[rd] &  &  & *+[F-]{6} \ar[rd] \ar[dd] & \\
*+[F-]{0} \ar[ru] \ar@(l,d)[] &  & *+[F-]{2} \ar@/^/[r] \ar[dl] & *+[F-]{1} \ar@/^/[l] \ar[ur] &  & *+[F-]{7} \ar[dl] \ar@(r,d)[] \\
 & *+[F-]{5} \ar[ul] \ar[uu] &  &  & *+[F-]{3} \ar[lll] \ar[ul] &
}
\end{align*}
\caption{The binary Collatz graphs $C(2)$ with modulus $4$ (left) and $C(3)$ with modulus $8$ (right).\label{fig2}}
\end{figure}

Upon further inspection, it turns out that in general, these graphs do not look particularly nice or structured. But when we take the modulus $m$ to be some power of $2$, these graphs do have a nice structure. From now on we will therefore focus on what we call \textit{binary modular Collatz graphs}, or simply \textit{binary Collatz graphs}. We write $C(k)$ for the modular Collatz graph with modulus $m = 2^k$, and we refer to this graph as the \textit{binary Collatz graph of dimension $k$}. For convenience, we write $T_k$ for forward iteration in the graph $C(k)$, i.e., $T_k(n)$ is the set of vertices $v$ in the graph $C(k)$ that are connected to $n$ by an edge $n \to v$. This relation can be explicitly written as \footnote{In this expression, numbers should be calculated modulo $2^k$, as all vertices correspond to congruence classes modulo $2^k$.}
\begin{align}
T_k(n) = \begin{cases} \{(3n + 1)/2, (3n + 1)/2 + 2^{k-1}\} & \text{if $n$ is odd,} \\ \{n/2, n/2 + 2^{k-1}\} & \text{if $n$ is even.} \end{cases}
\end{align}
For $k = 2, 3$ we get the graphs $C(2)$ and $C(3)$ shown in Figure~\ref{fig2}. These can also be found in a recent paper of Monks et al.~\cite[Figures~7.1, 7.2]{monks12}.

Looking at these graphs, we can immediately see a lot of structure. Both graphs have several symmetries, every vertex has two incoming and two outgoing edges, and reversing the direction of each edge leads to a graph isomorphic to the original graph. In the graph $C(2)$ we have also labeled each edge with the corresponding congruence class modulo $8$, e.g., the edge from $2$ to $3$ has a label $6$, because the numbers $a \equiv 2 \pmod 4$ satisfying $T(a) \equiv 3 \pmod 4$ are exactly all numbers $a \equiv 6 \pmod 8$. With this labeling, we can see a connection between $C(2)$ and $C(3)$: the latter can be formed by taking the so-called \textit{line graph} of the former, associating edges in $C(2)$ to vertices in $C(3)$ and connected edges in $C(2)$ to edges in $C(3)$.

\begin{figure}[t]
\centering
\begin{align*}
\xymatrix@C=1.3pc{
 & *+[F-]{001} \ar[rrr] \ar[rd] &  &  & *+[F-]{011} \ar[rd] \ar[dd] & \\
*+[F-]{000} \ar[ru] \ar@(l,d)[] &  & *+[F-]{010} \ar@/^/[r] \ar[dl] & *+[F-]{101} \ar@/^/[l] \ar[ur] &  & *+[F-]{111} \ar[dl] \ar@(r,d)[] \\
 & *+[F-]{100} \ar[ul] \ar[uu] &  &  & *+[F-]{110} \ar[lll] \ar[ul] &
} \quad \quad \quad \quad \quad \xymatrix{
 & *+[F-]{4} \ar[rrr] \ar[rd] &  &  & *+[F-]{6} \ar[rd] \ar[dd] & \\
*+[F-]{0} \ar[ru] \ar@(l,d)[] &  & *+[F-]{2} \ar@/^/[r] \ar[dl] & *+[F-]{5} \ar@/^/[l] \ar[ur] &  & *+[F-]{7} \ar[dl] \ar@(r,d)[] \\
 & *+[F-]{1} \ar[ul] \ar[uu] &  &  & *+[F-]{3} \ar[lll] \ar[ul] &
}
\end{align*}
\caption{The binary De Bruijn graph of dimension $3$, with binary (left) and numeric (right) labels on vertices.\label{fig3}}
\end{figure}
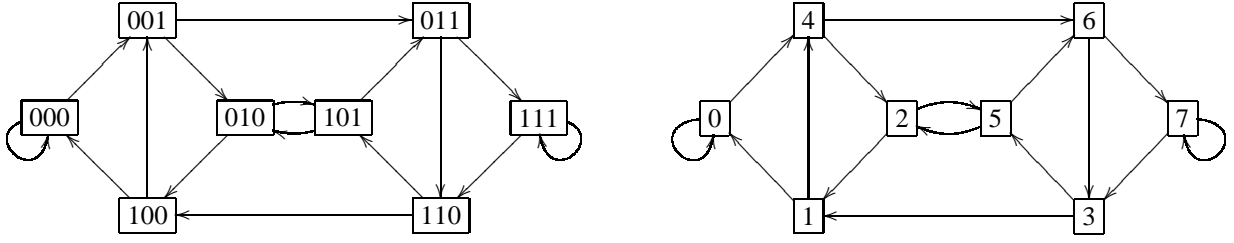

Seeing these beautiful graphs, one may wonder what is known about these graphs, and the title of this paper gives most of it away. Given an alphabet $\Sigma = \{0, 1, \ldots, p-1\}$ of size $p$ and a wordlength $k$, the \textit{$p$-ary De Bruijn graph of dimension $k$} \cite{debruijn46} is defined as the graph $B(p,k) = (V,E)$ with vertex set $V = \Sigma^k$, and an edge runs from the word $a_0 a_1 \cdots a_{k-1}$ to the word $b_0 b_1 \cdots b_{k-1}$ if and only if $a_{i+1} = b_i$ for $i = 0, \ldots, k - 2$. Thus, an edge runs from one word to another if the last $k - 1$ symbols of the first word overlap with the first $k-1$ symbols of the second word. When $p = 2$ we also refer to these graphs as \textit{binary De Bruijn graphs}. Besides viewing the vertices as words of a fixed length over some finite alphabet, it is also convenient to associate numbers between $0$ and $p^k - 1$ to the vertices. For this we identify words $b_0 b_1 \cdots b_{k-1}$ with numbers $\sum_{i=0}^{k-1} b_i p^i$. Figure~\ref{fig3} shows the two different labelings of the binary De Bruijn graph $B(2,3)$ of dimension $3$. For the remainder of this paper, we will choose to label the vertices with these numbers rather than with words over a finite alphabet.

We write $\sigma_{p,k}$ for forward iteration in the $p$-ary De Bruijn graph of dimension $k$. In terms of finite words over $\Sigma$, this corresponds to $\sigma_{p,k}(b_0 b_1 \cdots b_{k-1}) = \{b_1 b_2 \cdots b_{k-1} x: x \in \Sigma\}$. For now we will focus on the case $p = 2$, when the relation $\sigma_{2,k}$ can be described in terms of numbers as
\begin{align}
\sigma_{2,k}(n) = \begin{cases} \{(n - 1)/2, (n - 1)/2 + 2^{k-1}\} & \text{if $n$ is odd,} \\ \{n/2, n/2 + 2^{k-1}\} & \text{if $n$ is even.} \end{cases}
\end{align}
Since this relation shifts the bits to the left and appends a new bit, we will refer to this relation as the \textit{(binary) shift relation}.

De Bruijn graphs are closely related to, and mostly studied for, finding \textit{De Bruijn sequences}: cyclic sequences of symbols over a finite alphabet containing each sequence of length $k$ exactly once as a subsequence. For instance, a De Bruijn sequence for $k = 3$ is given by $00010111(00)$. These sequences correspond precisely to Hamiltonian paths in De Bruijn graphs. The above sequence corresponds to the path starting at $000$, and following the path $000 \stackrel{1}{\to} 001 \stackrel{0}{\to} 010 \stackrel{1}{\to} 101 \stackrel{1}{\to} 011 \stackrel{1}{\to} 111 \stackrel{(0)}{\to} 110 \stackrel{(0)}{\to} 100$, visiting each vertex of the graph exactly once, hence containing each sequence of length $3$ as a subsequence exactly once. Since De Bruijn graphs are Hamiltonian, such De Bruijn sequences exist for any value of $p$ and $k$, and can be constructed from a Hamiltonian path in a De Bruijn graph.\footnote{Due to the exponential size of the De Bruijn graphs, this method of finding De Bruijn sequences is not very efficient. In Section~\ref{sec:struct} we will see a more practical method for generating these sequences.}

Comparing Figures~\ref{fig2} and \ref{fig3}, it is clear that the graph $B(2,3)$ has exactly the same structure as the binary Collatz graph $C(3)$. In fact, for any value of $k$, the graphs $C(k)$ and $B(2,k)$ have exactly the same structure. This can be seen as follows. Let us use the notation $(x_i(n))$ for the $0-1$ sequence defined by $x_i(n) \equiv T^{i}(n) \pmod 2$, for $i \in \mathbb{N}$, as in \cite{lagarias85b}. For example, since the orbit of $n = 3$ under iterating $T$ is $(3, 5, 8, 4, 2, 1, 2, 1, 2, \ldots)$, the sequence $(x_i(n))$ is this orbit taken modulo $2$, i.e., $(1, 1, 0, 0, 0, 1, 0, 1, 0, \ldots)$. Then it is immediate that $x_i(T(n)) \equiv T^{i}(T(n)) = T^{i+1}(n) \equiv x_{i+1}(n) \pmod 2$, i.e., applying $T$ to $n$ means shifting the sequence $(x_i(n))$ one position to the left. The first $k$ terms of $(x_i(n))$ depend only on the value of $n \pmod{2^k}$, and $x_{k+1}(n)$ and $x_{k+1}(n + 2^k)$ always are different, as shown in \cite[Equation (2.9)]{lagarias85b}. So if we consider the (iterated) action of $T$ on a congruence class of some $n \pmod{2^k}$, then we find that the first $k$ elements of $(x_i(n))$ are identical for all numbers in this congruence class, and that the first $k$ elements of $(x_i(T(n)))$ are those of $(x_i(n))$ shifted by one to the left, with at the end added a $0$ or $1$, both occuring equally often. But this shows exactly that one forward iteration in the binary Collatz graph of dimension $k$ can be described as the shift relation on the bit sequences of length $k$, i.e., as the forward iteration on the binary De Bruijn graph of dimension $k$.

To summarize our result, we introduce the function $\Phi_k: \{0, \ldots, 2^k - 1\} \to \{0, \ldots, 2^k - 1\}$, defined by
\begin{align}
\Phi_k(n) = \sum_{i=0}^{k-1} x_i(n) 2^i.
\end{align}
Then we have the following statement.

\begin{theorem}
%
%
For any $k \geq 1$, the function $\Phi_k$
%
%
is an isomorphism between $C(k)$ and $B(2,k)$.
\end{theorem}

For the binary Collatz graphs, we introduced the notation $T_k$ to indicate forward iteration in the graph, and for binary De Bruijn graphs we wrote $\sigma_{2,k}$ for walking along edges in these graphs. Since $\Phi_k$ is an isomorphism between the two graphs, we can describe the relation between $T_k$ and $\sigma_{(2,k)}$ via $\Phi_k$ by
\begin{align}
T_k \equiv \Phi_k^{-1} \circ \sigma_{2,k} \circ \Phi_k,
\end{align}
where $\Phi_k^{-1}(\{a,b\}) = \{\Phi_k^{-1}(a), \Phi_k^{-1}(b)\}$. This is why $\Phi_k$ may be called a \emph{($k$-dimensional) conjugacy map}.

\begin{figure}[t]
\centering
\begin{align*}
\xymatrix@C=0.35cm@R=0.35cm{
 & & & & *+[F-]{12} \ar@/^/[rrrd] \ar@/^/[rdd] & & & & \\
 & *+[F-]{8} \ar[r] \ar@/^/[rrru] & *+[F-]{4} \ar[rr] \ar[dd] & & *+[F-]{2} \ar@/_/[dd] \ar[rr] & & *+[F-]{9} \ar[r] \ar[dl] & *+[F-]{14} \ar[dd] \ar[dr] & \\
*+[F-]{0} \ar[ru] \ar@(l,d)[] & & & *+[F-]{13} \ar[lu] \ar@/^/[uur] & & *+[F-]{6} \ar[dr] \ar@/^/[ldd] & & & *+[F-]{15} \ar[ld] \ar@(r,d)[] \\
 & *+[F-]{5} \ar[lu] \ar[uu] & *+[F-]{10} \ar[l] \ar[ru] & & *+[F-]{1} \ar[ll] \ar@/_/[uu] & & *+[F-]{11} \ar[ll] \ar[uu] & *+[F-]{7} \ar[l] \ar@/^/[llld] & \\
 & & & & *+[F-]{3} \ar@/^/[lllu] \ar@/^/[luu] & & & & \\
} \quad \quad \xymatrix@R=3pc{ & \\ \ar[r]^{\Phi_4} & } \quad \quad
\xymatrix@C=0.35cm@R=0.35cm{
 & & & & *+[F-]{12} \ar@/^/[rrrd] \ar@/^/[rdd] & & & & \\
 & *+[F-]{8} \ar[r] \ar@/^/[rrru] & *+[F-]{4} \ar[rr] \ar[dd] & & *+[F-]{10} \ar@/_/[dd] \ar[rr] & & *+[F-]{13} \ar[r] \ar[dl] & *+[F-]{14} \ar[dd] \ar[dr] & \\
*+[F-]{0} \ar[ru] \ar@(l,d)[] & & & *+[F-]{9} \ar[lu] \ar@/^/[uur] & & *+[F-]{6} \ar[dr] \ar@/^/[ldd] & & & *+[F-]{15} \ar[ld] \ar@(r,d)[] \\
 & *+[F-]{1} \ar[lu] \ar[uu] & *+[F-]{2} \ar[l] \ar[ru] & & *+[F-]{5} \ar[ll] \ar@/_/[uu] & & *+[F-]{11} \ar[ll] \ar[uu] & *+[F-]{7} \ar[l] \ar@/^/[llld] & \\
 & & & & *+[F-]{3} \ar@/^/[lllu] \ar@/^/[luu] & & & & \\
}
\end{align*}
\caption{The graphs $C(4)$ (left) and $B(2,4)$ (right), and the corresponding isomorphism $\Phi_4 = (1,5)(2,10)(9,13)$.\label{fig7}}
\end{figure}
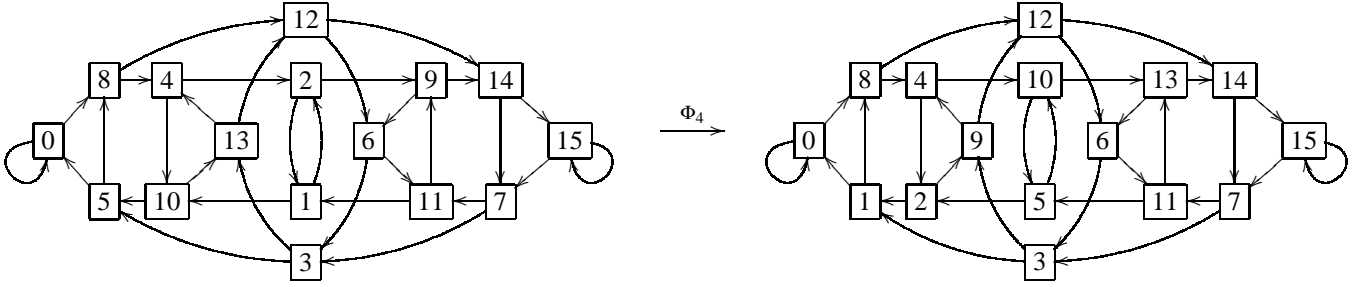

While the relation between De Bruijn graphs and the Collatz conjecture does not appear in the literature, these bijections $\Phi_k$ have been explicitly studied before by Bernstein and Lagarias \cite{bernstein94, bernstein96, lagarias85b}, but in a context different from graph isomorphisms. In \cite[Table 2]{lagarias85b} these bijections were explicitly given as permutations on numbers between $0$ and $2^k - 1$, e.g., for $k = 3$ we get the permutation $\Phi_3 \equiv (1,5)$ cf.\ Figures~\ref{fig2},~\ref{fig3}, and for $k = 4$ we get the permutation $\Phi_4 \equiv (1,5)(2,10)(9,13)$ cf.\ Figure~\ref{fig7}.

Since the graphs $C(k)$ and $B(2,k)$ are isomorphic, they share several properties. As mentioned before, the line graph of a $k$-dimensional binary Collatz graph is isomorphic to the $(k+1)$-dimensional binary Collatz graph, and the transpose graph of $C(k)$, obtained by reversing the direction of each edge, is again isomorphic to $C(k)$. One other remarkable property is that after exactly $k$ steps in these graphs, there is exactly one way to end up at any vertex in the graph (including the vertex itself). This can be expressed in terms of the adjacency matrix of the graph, $A_k$, by $(A_k^k)_{i,j} = 1$ for all indices $i,j$. This was previously noted by Feix et al.~\cite{feix94}, and Feix and Rouet~\cite{feix99}. Since every vertex has two outgoing edges, it follows that $(A_k^{\ell})_{i,j} = 2^{\ell - k}$ for $\ell \geq k$. In particular, this value does not depend on $i$ or $j$. Thus, if we start at any vertex in the graph, and take $k$ or more random steps in the graph, we can be anywhere with equal probability. So if we only know the $k$ least significant bits of a number $n$, then we know absolutely nothing about $T^{k}(n)$. This also gives some intuition why the $3n + 1$ problem is so hard: unless we know exactly what number we started out with, during each iteration we lose one bit of information about the resulting number.


We are ultimately looking for more insight into the $3n + 1$ problem, but these graphs are not quite the same as the Collatz graph. The vertices are congruence classes rather than numbers, and the graphs are finite. Still, there are relations between these binary Collatz graphs and the regular Collatz graph. For instance, if we simply restrict the Collatz graph to the vertex set $\{0, \ldots, 2^k - 1\}$, we get a subgraph of $C(k)$. Also, each edge in the Collatz graph corresponds to a sequence of edges in the graphs $C(k)$, e.g., the edge $5 \to 8$ in the Collatz graph corresponds to the edges $1 \to 0$ in $C(1)$ and $C(2)$, the edge $5 \to 0$ in $C(3)$, and the edge $5 \to 8$ in the graphs $C(k)$ for $k \geq 4$. By considering the sequence of graphs $C(k)$, and taking only those edges that appear in infinitely many of these graphs, we exactly get the Collatz graph on the natural numbers (including $0$).


\section{The $2$-adic Collatz graph and the $2$-adic De Bruijn graph}
\label{sec:infinite}

While considering finite Collatz graphs on congruence classes leads to nicely structured graphs, we would like to know more about the infinite Collatz graph on the natural numbers. Therefore it would be interesting to consider extensions of these binary Collatz graphs and binary De Bruijn graphs to larger (infinite) vertex sets, by in some way letting $k$ go to infinity. Then vertices may actually correspond to numbers rather than residue classes.

Although it is not exactly clear what ``$\lim_{k \to \infty} C(k)$'' means, we can naturally find an answer to this question by taking a detour along the De Bruijn graphs. First, we can easily extend the concept of De Bruijn graphs to infinite sequences of symbols over a finite alphabet. We define the \textit{infinite binary De Bruijn graph} by the graph $G = (V,E)$, where the vertex set is defined by $V = \{b_0 b_1 \cdots: b_i \in \{0,1\}\}$, and an edge runs from a vertex $a_0 a_1 \cdots$ to a vertex $b_0 b_1 \cdots$ if and only if $a_{i+1} = b_i$, for each $i \in \mathbb{N}$. Similar to the previous section, to these infinite sequences we may associate numbers, by considering the set of $2$-adic integers $\mathbb{Z}_2$ \cite{gouvea97}. We will identify sequences $b_0 b_1 \cdots$ with $2$-adic integers $\sum_{i=0}^{\infty} b_i 2^i$. With this labeling of the vertices, we will also refer to the infinite binary De Bruijn graph as the \textit{$2$-adic De Bruijn graph}, or $B(\mathbb{Z}_2)$. Note that in this graph, each vertex has only one outgoing edge. So while the shift relation $\sigma_{2,k}$ for finite binary De Bruijn graphs, mapping vertices $m$ to its set of neighbors, is not a proper function, forward iteration in the $2$-adic De Bruijn graph can be seen as a proper function from $\mathbb{Z}_2$ to $\mathbb{Z}_2$. We denote this function by $\sigma_2$, and it satisfies
\begin{align}
\sigma(n) = \begin{cases} (n - 1)/2 & \text{if $n$ is odd,} \\ n/2 & \text{if $n$ is even.}\end{cases}
\end{align}
In terms of binary sequences, this function $\sigma_2$ simply removes the first bit of a number and shifts the `remaining' bits one position to the left. This function is therefore known in the literature as the \textit{shift map} \cite{akin04, bernstein94, bernstein96, kraft10, monks04}.

Before going back to the $3n + 1$ problem, consider the limit of the isomorphisms $\Phi_k$, for $k \to \infty$. Using the definition $\Phi_k(n) = \sum_{i=0}^{k-1} x_i(n) 2^i$, we can just let $k$ tend to infinity to obtain $\Phi(n) = \sum_{i=0}^{\infty} x_i(n) 2^i$. In $\mathbb{Z}_2$, this is a convergent series, and for any $n \in \mathbb{Z}_2$, $\Phi(n)$ is a well-defined $2$-adic integer. By investigating $\Phi^{-1}(B(\mathbb{Z}_2))$ it then becomes clear what ``$\lim_{k \to \infty} C(k)$'' should be, as $\Phi^{-1}(B(\mathbb{Z}_2))$ corresponds exactly to the Collatz graph extended to the $2$-adic integers. The extension of the $3n + 1$ function $T$ to the $2$-adic integers was previously studied in \cite{akin04, bernstein94, bernstein96, kraft10, lagarias90, matthews92, monks04, wirsching98}. We will refer to this graph, with vertex set $\mathbb{Z}_2$ and edge set $\{n \to T(n): n \in \mathbb{Z}_2\}$, as the \textit{$2$-adic Collatz graph}, and denote it by $C(\mathbb{Z}_2)$. The following theorem is immediate.

\begin{theorem}
The graphs $C(\mathbb{Z}_2)$ and $B(\mathbb{Z}_2)$ are isomorphic, and the function $\Phi: \mathbb{Z}_2 \to \mathbb{Z}_2$, defined by
\begin{align}
\Phi(n) = \sum_{i=0}^{\infty} x_i(n) 2^i,
\end{align}
is an isomorphism from $C(\mathbb{Z}_2)$ to $B(\mathbb{Z}_2)$.
\end{theorem}

The function $\Phi$ is known in the literature as the \textit{conjugacy map} \cite{bernstein94, bernstein96, kraft10, lagarias85b, lagarias10, monks04, wirsching98}, and it clearly satisfies
\begin{align}
T \equiv \Phi^{-1} \circ \sigma \circ \Phi.
\end{align}
For values of $n$ for which the behaviour of iterates of $n$ is completely known, one can easily compute $\Phi(n)$. For instance, $\Phi(1) = 101010\ldots = -1/3$, $\Phi(5) = -13/3$, and all numbers ending in the cycle $\{1,2\}$ correspond to rational numbers in the $2$-adic De Bruijn graph with denominator $3$. The Collatz conjecture is equivalent to the statement $\Phi(\mathbb{N}_{+}) \subseteq \frac{1}{3}\mathbb{Z} \setminus \mathbb{Z}$ \cite{lagarias85b}.


\section{Structure inside the $2$-adic Collatz graph}
\label{sec:struct}

Let us now further investigate the graph $C(\mathbb{Z}_2)$. First, the graph $C(\mathbb{Z}_2)$ is not connected, and contains uncountably many components. It contains two types of components: cyclic components, corresponding to cycles of the function $T$; and divergent components, which do not contain a cycle, but extend infinitely far in both directions.

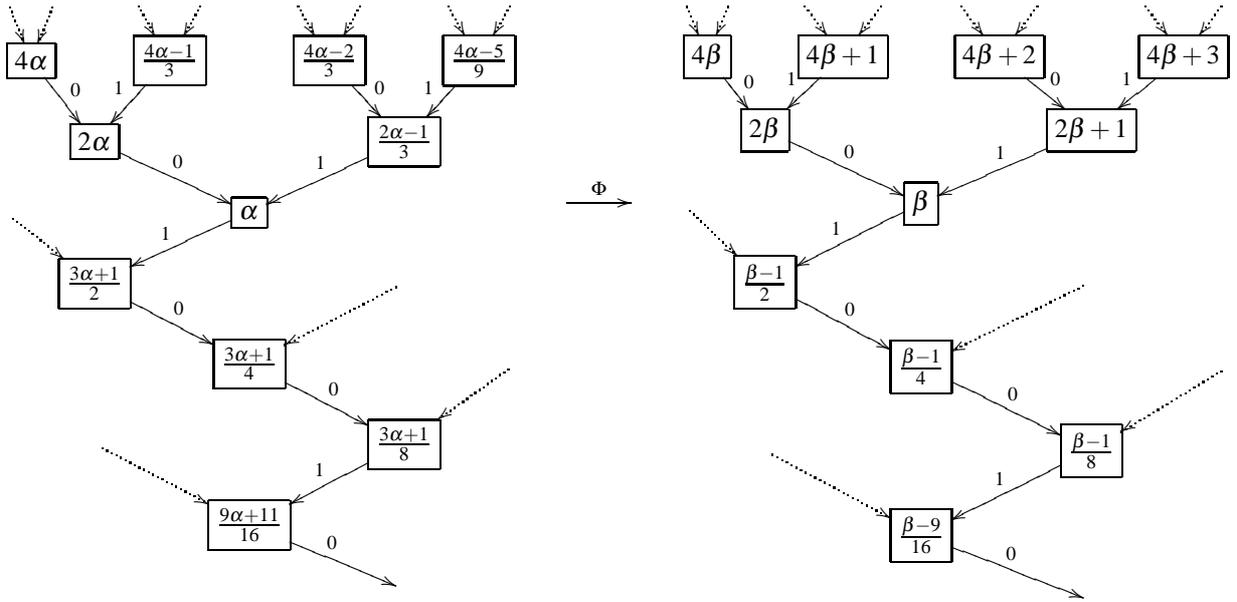
\begin{figure}[t]
\centering
\begin{align*}
\xymatrix@R=1pc@C=-.2pc{
 \ar@{.>}[rd] & & \ar@{.>}[ld] & & \ar@{.>}[rd] & & \ar@{.>}[ld] & & \ar@{.>}[rd] & & \ar@{.>}[ld] & & \ar@{.>}[rd] & & \ar@{.>}[ld] \\
 & *+[F-]{4\alpha} \ar[rrd]^{0} & & & & *+[F-]{\frac{4\alpha - 1}{3}} \ar[lld]_{1} & & & & *+[F-]{\frac{4\alpha - 2}{3}} \ar[rrd]^{0} & & & & *+[F-]{\frac{4\alpha - 5}{9}} \ar[lld]_{1} & \\
 & & & *+[F-]{2\alpha} \ar[rrrrd]^{0} & & & & & & & & *+[F-]{\frac{2\alpha - 1}{3}} \ar[lllld]_{1} & & & \\
\ar@{.>}[rrrd] & & & & & & & *+[F-]{\alpha} \ar[lllld]_{1} & & & & & & & \\
 & & & *+[F-]{\frac{3\alpha+1}{2}} \ar[rrrrd]^{0} & & & & & & & & \ar@{.>}[lllld] & & & \\
 & & & & & & & *+[F-]{\frac{3\alpha+1}{4}} \ar[rrrrd]^{0} & & & & & & & \ar@{.>}[llld] \\
 & & & \ar@{.>}[rrrrd] & & & & & & & & *+[F-]{\frac{3\alpha+1}{8}} \ar[lllld]_{1} & & & \\
 & & & & & & & *+[F-]{\frac{9\alpha+11}{16}} \ar[rrrrd]^{0} & & & & & & & \\
 & & & & & & & & & & & & & & \\
}
\quad \xymatrix@R=6pc{ & \\ \ar[r]^{\Phi} & } \quad
\xymatrix@R=1pc@C=-.2pc{
 \ar@{.>}[rd] & & \ar@{.>}[ld] & & \ar@{.>}[rd] & & \ar@{.>}[ld] & & \ar@{.>}[rd] & & \ar@{.>}[ld] & & \ar@{.>}[rd] & & \ar@{.>}[ld] \\
 & *+[F-]{4\beta} \ar[rrd]^{0} & & & & *+[F-]{4\beta + 1} \ar[lld]_{1} & & & & *+[F-]{4\beta + 2} \ar[rrd]^{0} & & & & *+[F-]{4\beta + 3} \ar[lld]_{1} & \\
 & & & *+[F-]{2\beta} \ar[rrrrd]^{0} & & & & & & & & *+[F-]{2\beta + 1} \ar[lllld]_{1} & & & \\
\ar@{.>}[rrrd] & & & & & & & *+[F-]{\beta} \ar[lllld]_{1} & & & & & & & \\
 & & & *+[F-]{\frac{\beta-1}{2}} \ar[rrrrd]^{0} & & & & & & & & \ar@{.>}[lllld] & & & \\
 & & & & & & & *+[F-]{\frac{\beta-1}{4}} \ar[rrrrd]^{0} & & & & & & & \ar@{.>}[llld] \\
 & & & \ar@{.>}[rrrrd] & & & & & & & & *+[F-]{\frac{\beta-1}{8}} \ar[lllld]_{1} & & & \\
 & & & & & & & *+[F-]{\frac{\beta-9}{16}} \ar[rrrrd]^{0} & & & & & & & \\
 & & & & & & & & & & & & & & \\
}
\end{align*}
\caption{A divergent component in $C(\mathbb{Z}_2)$ with an irrational number $\alpha = 10110\ldots$ (left), and the corresponding component in $B(\mathbb{Z}_2)$ with the irrational number $\Phi(\alpha) = \beta = 10010\ldots$ (right). Note that the first five bits of $\beta$ follow from $\Phi_5(10110) = 10010$, or equivalently, $\Phi_5(13) = 9$.\label{fig10}}
\end{figure}

The subgraph of $C(\mathbb{Z}_2)$ containing all divergent components, which we will denote by $C^{div}(\mathbb{Z}_2)$, is arguably the least interesting part of the graph. Although it contains uncountably many components, they are all pairwise isomorphic, and the vertices in each of these graphs are all indistinguishable, in the sense that the forward mapping is an automorphism. Figure~\ref{fig10} shows one of these components, containing an irrational $2$-adic integer $\alpha = 10110\cdots$, and the corresponding component in $B(\mathbb{Z}_2)$.

Note that theoretically, it is possible that a divergent component arises from some \textit{rational} $2$-adic integer, i.e., a rational number with an odd denominator. It has been conjectured \cite[Periodicity Conjecture]{lagarias85b} that all rational $2$-adic integers are part of a cyclic component, which can be formulated as $\Phi(\mathbb{Z}_2 \cap \mathbb{Q}) = \mathbb{Z}_2 \cap \mathbb{Q}$. But so far, we only know for certain that cyclic components must arise from rational numbers, i.e., $\Phi(\mathbb{Z}_2 \cap \mathbb{Q}) \supseteq \mathbb{Z}_2 \cap \mathbb{Q}$ \cite[Corollary 1]{bernstein94}.

\begin{figure}[t]
\centering
\begin{align*}
\xymatrix@C=1pc@R=1pc{
 & \ar@{.>}[dr] & & \ar@{.>}[dl] & \ar@{.>}[dr] & & \ar@{.>}[dl] & \\
\ar@{.>}[dr] & & *+[F-]{16} \ar[d] & & & *+[F-]{-\frac{11}{27}} \ar[d] & & \ar@{.>}[dl] \\
\ar@{.>}[r] & *+[F-]{5} \ar[r] & *+[F-]{8} \ar[d] & & & *+[F-]{-\frac{1}{9}} \ar[d] & *+[F-]{-\frac{2}{9}} \ar[l] & \ar@{.>}[l] \\
 & & *+[F-]{4} \ar[r] & *+[F-]{2} \ar@/^.5cm/[r] & *+[F-]{1} \ar@/^.5cm/[l] & *+[F-]{\frac{1}{3}} \ar[l] & & \\
\ar@{.>}[r] & *+[F-]{\frac{14}{3}} \ar[r] & *+[F-]{\frac{7}{3}} \ar[u] & & & *+[F-]{\frac{2}{3}} \ar[u] & *+[F-]{\frac{1}{9}} \ar[l] & \ar@{.>}[l] \\
\ar@{.>}[ur] & & *+[F-]{\frac{11}{9}} \ar[u] & & & *+[F-]{\frac{4}{3}} \ar[u] & & \ar@{.>}[ul] \\
 & \ar@{.>}[ur] & & \ar@{.>}[ul] & \ar@{.>}[ur] & & \ar@{.>}[ul] & \\
}
\quad \xymatrix@R=6pc{ & \\ \ar[r]^{\Phi} & } \quad
\xymatrix@C=1pc@R=1pc{
 & \ar@{.>}[dr] & & \ar@{.>}[dl] & \ar@{.>}[dr] & & \ar@{.>}[dl] & \\
\ar@{.>}[dr] & & *+[F-]{-\frac{16}{3}} \ar[d] & & & *+[F-]{\frac{13}{3}} \ar[d] & & \ar@{.>}[dl] \\
\ar@{.>}[r] & *+[F-]{-\frac{13}{3}} \ar[r] & *+[F-]{-\frac{8}{3}} \ar[d] & & & *+[F-]{\frac{5}{3}} \ar[d] & *+[F-]{\frac{10}{3}} \ar[l] & \ar@{.>}[l] \\
 & & *+[F-]{-\frac{4}{3}} \ar[r] & *+[F-]{-\frac{2}{3}} \ar@/^.5cm/[r] & *+[F-]{-\frac{1}{3}} \ar@/^.5cm/[l] & *+[F-]{\frac{1}{3}} \ar[l] & & \\
\ar@{.>}[r] & *+[F-]{-\frac{10}{3}} \ar[r] & *+[F-]{-\frac{5}{3}} \ar[u] & & & *+[F-]{\frac{2}{3}} \ar[u] & *+[F-]{\frac{7}{3}} \ar[l] & \ar@{.>}[l] \\
\ar@{.>}[ur] & & *+[F-]{-\frac{7}{3}} \ar[u] & & & *+[F-]{\frac{4}{3}} \ar[u] & & \ar@{.>}[ul] \\
 & \ar@{.>}[ur] & & \ar@{.>}[ul] & \ar@{.>}[ur] & & \ar@{.>}[ul] & \\
}
\end{align*}
\caption{The components corresponding to the cycle $(01)$, in $C(\mathbb{Z}_2)$ (left) and in $B(\mathbb{Z}_2)$ (right). The Collatz conjecture states that the component on the left contains all positive integers.\label{fig11}}
\end{figure}
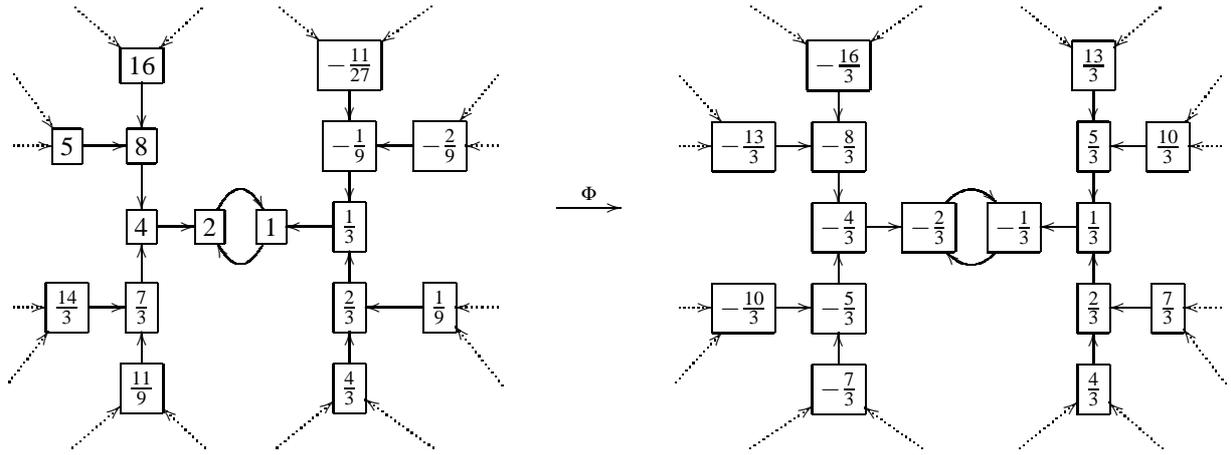

The other part of the graph, formed by the $2$-adic integers that end in a cycle, has different properties. The number of components is only countably infinite, since we can enumerate all possible cycles to find all components. For instance, there are two components with a cycle of length $1$ (containing the cycles $0 \stackrel{0}{\to} 0$ and $-1 \stackrel{1}{\to} -1$), and there is exactly one component with a cycle of length $2$, shown in Figure~\ref{fig11}. The Collatz conjecture states that this latter component contains the entire Collatz graph on the positive integers as a subgraph.

Writing $M_k$ for the number of components in $C(\mathbb{Z}_2)$ with a cycle of length $k$, by counting words of length $k$ and computing their cyclic orders, it follows that
\begin{align}
2^k = \sum_{d|k} d \ M_d.
\end{align}
Applying M\"{o}bius inversion to the above, we obtain a direct formula for $M_k$ as
\begin{align}
M_k = \frac{1}{k} \sum_{d|k} \mu(d) 2^{k/d} = \frac{1}{n} 2^n + O(2^{n/2}).
\end{align}
The sequence $(M_k)_{k \geq 1} = (1, 2, 1, 2, 3, 6, \ldots) $ has been encountered before in the context of the Collatz conjecture \cite{lagarias90}, but has also been studied in different contexts \cite{oeis12} and is known as Moreau's necklace-counting function \cite{moreau72}. What this function $M$ is actually counting is the number of Lyndon words of length $k$: strings of length $k$ that are inequivalent modulo cyclic rotations, and with period equal to $k$ \cite{lyndon54}. Somewhat surprisingly, a connection between Lyndon words and De Bruijn sequences has been made before in the literature \cite{fredricksen78, fredricksen86}: one may obtain a De Bruijn sequence of order $k$ by appropriately concatenating Lyndon words of length some divisor $d$ of $k$. For instance, for $k = 4$, the Lyndon words of lengths $d|k$ are given by $0, 1, 01, 0001, 0011, 0111$. By extending these words to words of length $4$, we may order them lexicographically as $0(000), 0001, 0011, 01(01), 0111, 1(111)$. Concatenating the Lyndon words in this order, we get the sequence
\begin{align}
0 | 0001 | 0011 | 01 | 0111 | 1,
\end{align}
which is indeed a De Bruijn sequence of order $4$, corresponding to a Hamiltonian path in the graph $B(2,4)$. This algorithm is known in the literature as the FKM algorithm, after its authors Fredricksen, Kessler and Maiorana \cite{fredricksen78, fredricksen86}.

Another fascinating property of the cyclic part of the $2$-adic Collatz graph, is the fact that there is a one-to-one correspondence between cycles in $C(\mathbb{Z}_2)$, and integer cycles of the family of $3n + b$ functions $T^{(3,b)}$ defined by
\begin{align}
T^{(3,b)}(n) = \begin{cases} (3n + b)/2 & \text{if $n$ is odd,} \\ n/2 & \text{if $n$ is even,} \end{cases}
\end{align}
for odd $b$ coprime to $3$. More precisely, if there is a rational cycle in $C(\mathbb{Z}_2)$ with all numbers having denominator $b$, then this cycle corresponds exactly to an integer cycle of the $3n + b$ problem. This follows from the fact that $T(n/b) = T^{(3,b)}(n) / b$. This correspondence was previously noted by Lagarias \cite{lagarias90}. So together, the $3n + b$ problems on integers are very structured, as we know that every cycle (Lyndon word) corresponds to exactly one value $b$, and finding the $b$ associated to a given Lyndon word is also not so hard. But when we try to find all cycles associated to a value $b$, we get stuck. Note that solving this problem would solve the Collatz conjecture, and much more.

Going back to the complete $2$-adic Collatz graph $C(\mathbb{Z}_2)$, we saw that we may divide the graph in two parts: a cyclic part, and a divergent part. Furthermore, we know that all cycles must correspond to rational numbers, but we do not know whether all rational numbers also correspond to cycles. So theoretically, we have three different types of components: cyclic, rational components; divergent, rational components; and divergent, irrational components. We will denote these parts of the graph with $C^{cyc}(\mathbb{Z}_2 \cap \mathbb{Q})$, $C^{div}(\mathbb{Z}_2 \cap \mathbb{Q})$, and $C(\mathbb{Z}_2 \setminus \mathbb{Q})$ respectively. Each of these parts of the graph corresponds to a part of $B(\mathbb{Z}_2)$, under the map $\Phi$.

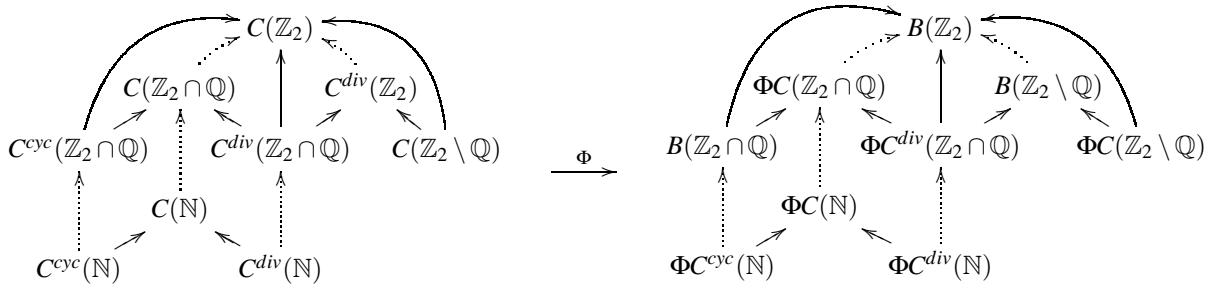
\begin{figure}[t]
\centering
\begin{align*}
\xymatrix@C=-1.3pc@R=.5pc{
 &  & C(\mathbb{Z}_2) & & \\
 & C(\mathbb{Z}_2 \cap \mathbb{Q}) \ar@{.>}@/^/[ur] & & C^{div}(\mathbb{Z}_2) \ar@{.>}@/_/[ul] & \\
C^{cyc}(\mathbb{Z}_2 \cap \mathbb{Q}) \ar[ru] \ar@/^2pc/[rruu] & & C^{div}(\mathbb{Z}_2 \cap \mathbb{Q}) \ar[lu] \ar[ru] \ar[uu] & & C(\mathbb{Z}_2 \setminus \mathbb{Q}) \ar[lu] \ar@/_2pc/[lluu] \\
 & C(\mathbb{N}) \ar@{.>}[uu] & & & \\
C^{cyc}(\mathbb{N}) \ar[ru] \ar@{.>}[uu] & & C^{div}(\mathbb{N}) \ar[lu] \ar@{.>}[uu] & & &
}
\quad \xymatrix@R=4pc{ & \\ \ar[r]^{\Phi} & } \quad
\xymatrix@C=-1.3pc@R=.5pc{
 & & B(\mathbb{Z}_2) & & \\
 & \Phi C(\mathbb{Z}_2 \cap \mathbb{Q}) \ar@{.>}@/^/[ur] & & B(\mathbb{Z}_2 \setminus \mathbb{Q}) \ar@{.>}@/_/[ul] & \\
B(\mathbb{Z}_2 \cap \mathbb{Q}) \ar[ru] \ar@/^2.5pc/[rruu] & & \Phi C^{div}(\mathbb{Z}_2 \cap \mathbb{Q}) \ar[lu] \ar[ru] \ar[uu] & & \Phi C(\mathbb{Z}_2 \setminus \mathbb{Q}) \ar[lu] \ar@/_2pc/[lluu] \\
 & \Phi C(\mathbb{N}) \ar@{.>}[uu] & & & \\
\Phi C^{cyc}(\mathbb{N}) \ar[ru] \ar@{.>}[uu] & & \Phi C^{div}(\mathbb{N}) \ar[lu] \ar@{.>}[uu] & &
}
\end{align*}
\caption{The structure inside $C(\mathbb{Z}_2)$ (left) and the corresponding structure inside $B(\mathbb{Z}_2)$ (right). Solid arrows pointing to the same larger graph imply that the large graph can be partitioned into these smaller graphs. Dotted lines correspond to graph inclusions.\label{fig12}}
\end{figure}

The complete structure of $C(\mathbb{Z}_2)$ and its image under $\Phi$ in $B(\mathbb{Z}_2)$ are summarized in Figure~\ref{fig12}. The periodicity conjecture \cite{lagarias85b} states that the middle part, $C^{div}(\mathbb{Z}_2 \cap \mathbb{Q})$, is empty. If this is true, then the picture somewhat simplifies, as then we would get $\Phi C(\mathbb{Z}_2 \cap \mathbb{Q}) = B(\mathbb{Z}_2 \cap \mathbb{Q})$ and $\Phi C(\mathbb{Z}_2 \setminus \mathbb{Q}) = B(\mathbb{Z}_2 \setminus \mathbb{Q})$.

The challenge of the Collatz conjecture now can be seen as a study in more detail of the embedding of $C(\mathbb{N}_+)$ into $C(\mathbb{Z}_2)$. The idea is that $C(\mathbb{N}_+)$ is a rather chaotic object, which lies somehow inside the beautifully structured object $C(\mathbb{Z}_2)$. The latter graph is isomorphic via the conjugacy map $\Phi$ to $B(\mathbb{Z}_2)$, which indeed has a rich and well understood structure. Similarly other interesting subgraphs of $C(\mathbb{Z}_2)$ can via $\Phi$ be embedded into $B(\mathbb{Z}_2)$.  


\section{Generalizations and the $p$-adic De Bruijn graphs}
\label{sec:generalizations}

\subsection{$an + b$ problems and the $2$-adic De Bruijn graph}
\label{sub:an+b}

In the previous section we briefly mentioned a generalization of the $3n + 1$ function $T$ to $3n + b$ functions $T^{(3,b)}$. A frequently considered further generalization of these functions is the family of \textit{$an + b$ functions} \cite{akin04, bernstein94, bernstein96, kraft10, lagarias90, matthews92, monks04, wirsching98}. For odd $a,b \in \mathbb{Z}$, we define the $an + b$ function $T^{(a,b)}$ by
\begin{align}
T^{(a,b)}(n) = \begin{cases} (an + b)/2 & \text{if $n$ is odd,} \\ n/2 & \text{if $n$ is even}. \end{cases} \label{eq:gen1}
\end{align}
For $a = 3$, we expect that no divergent paths exist, since with equal `probability', we either multiply a number by $1/2$ or (roughly) by $3/2$. So on average, the numbers increase by a factor $\sqrt{3/4} < 1$. For $a \geq 5$, this heuristic suggests a different behaviour. For instance, for the $5n + 1$ problem and the $7n + 1$ problem, we expect most paths to diverge \cite{crandall78, lagarias10}.

Although in this sense, the $an + b$ problems behave differently from the $3n + 1$ problem, the story described above related to De Bruijn graphs applies to all of these problems as well. Defining $C^{(a,b)}(\cdot)$ similarly as the Collatz graphs $C(\cdot)$ but for the $an + b$ function $T^{(a,b)}$, and defining $x_i^{(a,b)}(n)$ analogously to $x_i(n)$, we obtain the following result.
\begin{theorem}
For any $k \geq 1$, the graphs $C^{(a,b)}(k)$ and $B(2,k)$ are isomorphic, and the function $\Phi^{(a,b)}_k: \{0, \ldots, 2^k - 1\} \to \{0, \ldots, 2^k - 1\}$, defined by
\begin{align}
\Phi^{(a,b)}_k(n) = \sum_{i=0}^{k-1} x^{(a,b)}_i(n) 2^i,
\end{align}
is an isomorphism between $C^{(a,b)}(k)$ and $B(2,k)$. Furthermore, the graphs $C^{(a,b)}(\mathbb{Z}_2)$ and $B(\mathbb{Z}_2)$ are isomorphic, and the $an + b$ conjugacy map $\Phi^{(a,b)}: \mathbb{Z}_2 \to \mathbb{Z}_2$, defined by
\begin{align}
\Phi^{(a,b)}(n) = \sum_{i=0}^{\infty} x^{(a,b)}_i(n) 2^i,
\end{align}
is an isomorphism from $C^{(a,b)}(\mathbb{Z}_2)$ to $B(\mathbb{Z}_2)$, satisfying
\begin{align}
T^{(a,b)} \equiv (\Phi^{(a,b)})^{-1} \circ \sigma_2 \circ \Phi^{(a,b)}.
\end{align}
\end{theorem}

For example, for the $5n + 1$ problem and $k = 3$, we get the isomorphism $\Phi_3^{(5,1)} \equiv (1,3)(2,6)(5,7)$.

It is of interest to observe that for any $an + b$ problem, the binary modular graphs that appear are always isomorphic to the same binary De Bruijn graphs. Only the labeling of the graphs depends on the choice of $a$ and $b$. In other words, the binary and $2$-adic De Bruijn graphs themselves do not contain much information anymore about the generalized Collatz problems. The wanted information is contained in the isomorphisms (the labelings of the graph edges, the conjugacy maps), which are different for each $a, b$, and not well understood at all. That these maps really are different can be seen by noticing that the conjugacy map for the $3n + 1$ function maps $C(\mathbb{N}_+)$ into (most probably) one cyclic connected component in $B(\mathbb{Z}_2)$, while for the $5n + 1$ function the image of $C^{(5,1)}(\mathbb{N}_+)$ under its conjugacy map ends up in at least 3 cyclic components, but (most probably) in infinitely many components, most of which are divergent. While $C^{div}(\mathbb{Z}_2\cap\mathbb{Q})$ is conjectured to be empty, $(C^{(5,1)})^{div}(\mathbb{Z}_2\cap\mathbb{Q})$ is conjectured to have infinitely many distinct components. This may be seen as an illustration of the profound difficulty of the Collatz conjecture and its siblings.

Note that for $a = 1$ and $b = -1$, the isomorphism $\Phi^{(1,-1)}$ is extremely well understood, since $T^{(1,-1)} \equiv \sigma_2$ and $\Phi^{(1,-1)} \equiv \mathrm{id}$, the identity map. So one could also write $B(2,k) = C^{(1,-1)}(k)$ and $B(\mathbb{Z}_2) = C^{(1,-1)}(\mathbb{Z}_2)$.


\subsection{$p$-ary functions and the $p$-adic De Bruijn graphs}
\label{sub:p}

\begin{figure}[t]
\centering
\begin{align*}
\xymatrix@C=2pc@R=1.2pc{
 & *+[F-]{1} \ar[dddl] \ar@(l,u)[] \ar[ddrr] & & & *+[F-]{7} \ar[lll] \ar@/^/[ddl] \ar[llllddd] & \\
 & & & & & \\
 & & & *+[F-]{6} \ar@/^/[uur] \ar[rrd] \ar[dd] & & \\
*+[F-]{5} \ar[rddd] \ar@/^/[rr] \ar[rrrrddd] & & *+[F-]{4} \ar@/^/[ll] \ar[uuul] \ar[ur] & & & *+[F-]{0} \ar@(r,d)[] \ar[uuul] \ar[lld] \\
 & & & *+[F-]{2} \ar[lu] \ar[lldd] \ar@/^/[rdd] & & \\
 & & & & & \\
 & *+[F-]{8} \ar[rrr] \ar[ruuu] \ar@(l,d)[] & & & *+[F-]{3} \ar[uuur] \ar[uuuuuu] \ar@/^/[uul] &
}
\quad \quad \xymatrix{ & \\ \ar[r]^{\Phi_{3,2}^{(f_0)}} & } \quad \quad
\xymatrix@C=2pc@R=1.2pc{
 & *+[F-]{4} \ar[dddl] \ar@(l,u)[] \ar[ddrr] & & & *+[F-]{3} \ar[lll] \ar@/^/[ddl] \ar[llllddd] & \\
 & & & & & \\
 & & & *+[F-]{1} \ar@/^/[uur] \ar[rrd] \ar[dd] & & \\
*+[F-]{7} \ar[rddd] \ar@/^/[rr] \ar[rrrrddd] & & *+[F-]{5} \ar@/^/[ll] \ar[uuul] \ar[ur] & & & *+[F-]{0} \ar@(r,d)[] \ar[uuul] \ar[lld] \\
 & & & *+[F-]{6} \ar[lu] \ar[lldd] \ar@/^/[rdd] & & \\
 & & & & & \\
 & *+[F-]{8} \ar[rrr] \ar[ruuu] \ar@(l,d)[] & & & *+[F-]{2} \ar[uuur] \ar[uuuuuu] \ar@/^/[uul] &
}
\end{align*}
\caption{The ternary modular graph of dimension $2$ corresponding to the function $f_0$ in \eqref{eq:gen3} and the ternary De Bruijn graph of dimension $2$, $B(3,2)$. The corresponding isomorphism $\Phi^{(f_0)}_{3,2}$ of order $7$ is given by the permutation $\Phi^{(f_0)}_{3,2} \equiv (1,4,5,7,3,2,6)$. \label{fig16}}
\end{figure}
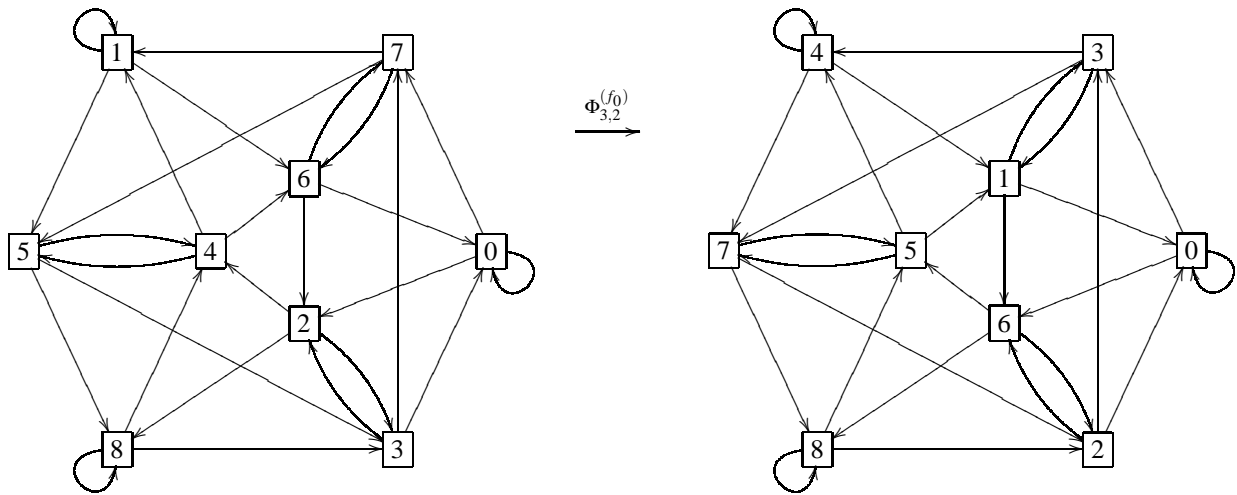

A further generalization of the function $T$ was considered in, e.g., \cite{feix94, feix99, matthews84, matthews85, matthews92, wirsching98}. For some appropriately chosen integers $a_i$ and $b_i$, let $f$ be defined by
\begin{align}
f(n) = \begin{cases} (a_0 n + b_0)/p & \text{if $n \equiv 0 \pmod p$,} \\ (a_1 n + b_1)/p & \text{if $n \equiv 1 \pmod p$,} \\ \ldots & \ldots \\ (a_{p-1} n + b_{p-1})/p & \text{if $n \equiv p-1 \pmod p$.} \end{cases} \label{eq:gen2}
\end{align}
For instance, a problem Lagarias attributes to Collatz \cite{lagarias85b} concerns the function $f_0$ defined by
\begin{align}
f_0(n) = \begin{cases} 2n/3 & \text{if $n \equiv 0 \pmod 3$,} \\ (4n-1)/3 & \text{if $n \equiv 1 \pmod 3$,} \\ (4n+1)/3 & \text{if $n \equiv 2 \pmod 3$.} \end{cases} \label{eq:gen3}
\end{align}
Iterating these functions leads to similar behaviour as the $an + b$ problems. In this case, when considering graphs on congruence classes modulo $p^k$, we find a relation with $p$-ary and $p$-adic De Bruijn graphs as follows.
\begin{theorem}
For any $k \geq 1$, the graphs $C^{(f)}(p,k)$ and $B(p,k)$ are isomorphic, and the function $\Phi^{(f)}_{p,k}: \{0, \ldots, p^k - 1\} \to \{0, \ldots, p^k - 1\}$, defined by
\begin{align}
\Phi^{(f)}_{p,k}(n) = \sum_{i=0}^{k-1} x^{(f)}_i(n) p^i,
\end{align}
is an isomorphism between $C^{(f)}(p,k)$ and $B(p,k)$. Furthermore, the graphs $C^{(f)}(\mathbb{Z}_p)$ and $B(\mathbb{Z}_p)$ are isomorphic, and the function $\Phi_p^{(f)}: \mathbb{Z}_p \to \mathbb{Z}_p$, defined by
\begin{align}
\Phi_p^{(f)}(n) = \sum_{i=0}^{\infty} x^{(f)}_i(n) p^i,
\end{align}
is an isomorphism from $C^{(f)}(\mathbb{Z}_p)$ to $B(\mathbb{Z}_p)$, satisfying
\begin{align}
f \equiv (\Phi_p^{(f)})^{-1} \circ \sigma_p \circ \Phi_p^{(f)}.
\end{align}
\end{theorem}
As an example, Figure~\ref{fig16} shows the $2$-dimensional ternary graph corresponding to the function $f_0$ in Equation~\eqref{eq:gen3}, and the corresponding ternary De Bruijn graph.

\end{document}